\documentclass[12pt,leqno]{article}
\usepackage{amsfonts}
\usepackage{amsmath, amsthm, amsfonts, amssymb, color}
\usepackage{mathrsfs}
\usepackage{color}
\newtheorem{thm}{Theorem}[section]

\theoremstyle{definition}

\newcommand{\scr}[1]{\mathscr #1}
\definecolor{wco}{rgb}{0.5,0.2,0.3}

\numberwithin{equation}{section} \theoremstyle{remark}

\newcommand{\ua}{\uparrow}

\title{{\bf Asymptotics in   Wasserstein Distance  for  Empirical Measures of Markov Processes}\footnote{Supported in
 part by  NNSFC (11921001) and the National Key R\&D Program of China (No. 2020YFA0712900).} }
\author{{\bf  Feng-Yu Wang   }\\
\footnotesize{ Center for Applied Mathematics, Tianjin University, Tianjin 300072, China}
\footnotesize{    wangfy@tju.edu.cn}}
\begin{document}
\allowdisplaybreaks
\def\R{\mathbb R}  \def\ff{\frac} \def\ss{\sqrt} \def\B{\mathbf
B}
\def\N{\mathbb N} \def\kk{\kappa} \def\m{{\bf m}}
\def\ee{\varepsilon}\def\ddd{D^*}
\def\dd{\delta} \def\DD{\Delta} \def\vv{\varepsilon} \def\rr{\rho}
\def\<{\langle} \def\>{\rangle}
  \def\nn{\nabla} \def\pp{\partial} \def\E{\mathbb E}
\def\d{\text{\rm{d}}} \def\bb{\beta} \def\aa{\alpha} \def\D{\scr D}
  \def\si{\sigma} \def\ess{\text{\rm{ess}}}\def\s{{\bf s}}
\def\beg{\begin} \def\beq{\begin{equation}}  \def\F{\scr F}
\def\Ric{\mathcal Ric} \def\Hess{\text{\rm{Hess}}}
\def\e{\text{\rm{e}}} \def\ua{\underline a} \def\OO{\Omega}  \def\oo{\omega}
 \def\tt{\tilde}\def\[{\lfloor} \def\]{\rfloor}
\def\cut{\text{\rm{cut}}} \def\P{\mathbb P} \def\ifn{I_n(f^{\bigotimes n})}
\def\C{\scr C}      \def\aaa{\mathbf{r}}     \def\r{r}
\def\gap{\text{\rm{gap}}} \def\prr{\pi_{{\bf m},\nu}}  \def\r{\mathbf r}
\def\Z{\mathbb Z} \def\vrr{\nu} \def\ll{\lambda}
\def\L{\scr L}\def\Tt{\tt} \def\TT{\tt}\def\II{\mathbb I}
\def\i{{\rm in}}\def\Sect{{\rm Sect}}  \def\H{\mathbb H}
\def\M{\mathbb M}\def\Q{\mathbb Q} \def\texto{\text{o}} \def\LL{\Lambda}
\def\Rank{{\rm Rank}} \def\B{\scr B} \def\i{{\rm i}} \def\HR{\hat{\R}^d}
\def\to{\rightarrow} \def\gg{\gamma}
\def\EE{\scr E} \def\W{\mathbb W}
\def\A{\scr A} \def\Lip{{\rm Lip}}\def\S{\mathbb S}
\def\BB{\scr B}\def\Ent{{\rm Ent}} \def\i{{\rm i}}\def\itparallel{{\it\parallel}}
\def\g{{\mathbf g}}\def\Sect{{\mathcal Sec}}\def\T{\mathcal T}\def\BB{{\bf B}}
\def\f\ell \def\g{\mathbf g}\def\BL{{\bf L}}  \def\BG{{\mathbb G}}
\def\Bd{{D^E}} \def\BdP{D^E_\phi} \def\Bdd{{\bf \dd}} \def\Bs{{\bf s}} \def\GA{\scr A}
\def\Bg{{\bf g}}  \def\Bdd{\psi_B} \def\supp{{\rm supp}}\def\div{{\rm div}}
\def\ddiv{{\rm div}}\def\osc{{\bf osc}}\def\1{{\bf 1}}\def\BD{\mathbb D}
\def\H{{\bf H}}\def\gg{\gamma} \def\n{{\mathbf n}} \def\GG{\Gamma}
\maketitle

\begin{abstract} In this paper we introduce some recent progresses on the convergence rate in Wasserstein distance for   empirical measures of Markov processes.
For diffusion processes on compact manifolds possibly with reflecting or killing boundary conditions, the sharp convergence rate as well as renormalization limits are presented in terms of the dimension
of the manifold and the spectrum of the generator.  For general ergodic Markov processes, explicit  estimates   are presented for the convergence rate by using a nice reference diffusion process,
which are illustrated by some typical examples.  Finally, some techniques are introduced to estimate the Wasserstein distance between empirical measures.  \end{abstract} \noindent

 \vskip 2cm

 \section{Introduction}

The empirical measure is a fundamental statistic to estimate the stationary distribution of a Markov process.   In this paper, we study the long time behavior of empirical measures for Markov processes under
the Wasserstein distance. For two nonnegative  functions $f,g$  on a space $E$, we write $f\lesssim g$ if there exists a constant $c>0$ such that $f\le c g$ holds on $E$, and write $f\asymp g$ if $f\lesssim g$ and $g\lesssim f$.

Let $(M,\rr)$ be a Polish space, let $\scr P$ be the space of all probability measures on $M$.
For a Markov process  $X_t$  on $M$, the    empirical measure  is defined as
$$\mu_t:=\ff 1 t \int_0^t \dd_{X_s}\d s$$ for $t>0$, where $\dd_{X_s}$ is the Dirac measure at $X_s$.
We intend to study the long time behavior of $\mu_t$ under the $p$-Wasserstein distance for any $p\in [1,\infty)$:
$$\W_p(\mu_1,\mu_2):= \inf_{\pi\in \scr C(\mu_1,\mu_2)}\bigg( \int_{M\times M} \rr^p\d\pi\bigg)^{\ff 1 p},\ \ \mu_1,\mu_2\in \scr P,$$
where $\scr C(\mu_1,\mu_2)$ is the set of all couplings of $\mu_1$ and $\mu_2$.

To this end, we will  consider the ergodicity   and quasi-ergodicity cases respectively.

\subsection{Ergodicity case}

For any $\nu\in \scr P$, let $\E^\nu$ denote the expectation for the Markov process $X_t$ with initial distribution $\nu$, and denote by  $P_t^*\nu$   the law of $X_t$ with initial distribution $\nu$. We have
$$\int_M f\d(P_t^*\nu)= \E^\nu[f(X_t)] =\int_M P_tf(x) \nu(\d x),\ \ \ f\in \B_b(M)$$
where $P_t f(x):= \E^x [f(X_t)]$, and $\B_b(M)$ the   class of all bounded measurable functions on $M$.

If  $\mu\in \scr P$ satisfies $P_t^*\mu=\mu$ for all $t\ge 0$, we call $\mu$ an invariant probability measure of the Markov process. If furthermore
$$\lim_{t\to\infty} P_t^*\nu=\mu\ \text{weakly} $$
holds for any $\nu\in \scr P$, we call the Markov process ergodic.

In particular, when the Markov process is exponentially ergodic in $L^2(\mu)$, i.e.
$$\|P_t-\mu\|_{L^2(\mu)}\le c\e^{-\ll t},\ \ t\ge 0$$ holds for some constant $c,\ll>0$, where $\|\cdot\|_{L^2(\mu)}$ is the operator norm in $L^2(\mu)$ and
$  \mu(f):=\int_Ef\d\mu$  for $f\in L^2(\mu)$, we have
$${\bf V}_f:=\int_0^\infty \mu(\hat fP_t  \hat f)\d t\in (0,\infty),\ \  0\ne \hat f:= f-\mu(f), \ f \in L^2(\mu).$$
Moreover, according to \cite{Wu2}, for any  $f\in L^2(\mu)$, $\P$-a.s.
$$\lim_{t\to\infty} \mu_t(f)=\mu(f),$$
and
$$\lim_{t\to\infty}\ss t \big(\mu_t(f)-\mu(f)\big)=N(0, 2{\bf V}_f)\ \ \text{in\ law},$$
where $N(0,2{\bf V}_f)$ is the normal distribution with mean $0$ and variance $2{\bf V}_f.$

Noting that for any $p\ge 1$
$$\W_p(\mu_t,\mu)\ge \W_1(\mu_t,\mu)=\sup_{\|f\|_{Lip} \le 1} |\mu_t(f)-\mu(f)|,$$
where $\|\cdot\|_{Lip}$ is the Lipschitz constant,
the above result implies that the convergence of $\E [\W_p(\mu_t,\mu)]$ could not faster than $t^{-\ff 1 2}$ as $t\to\infty$.

We will show that for elliptic diffusions on compact manifolds, this best possible rate is reached if and only if the dimension of manifold is not larger than $3$, and in this case we also derive the exact   limit for
$t\E[\W_2(\mu_t,\mu)^2]$ as $t\to\infty$, which is given explicitly by the spectrum of the generator.  

\subsection{Quasi-ergodicity case }

Consider Markov process $X_t$ with finite life time $\tau$. Typical examples are the diffusion processes on a domain with killing boundary, for which the life time is the first hitting time to the boundary.
In this case, we study the long time behavior of the empirical measure $\mu_t$ under the condition $\{t<\tau\}$, i.e.   the process   is not yet killed at time $t$.

A probability measure $\mu$ is called quasi-invariant for $X_t$, if
$$\E^\mu(f(X_t)|t<\tau)= \mu(f), \ \ t\ge 0, f\in \B_b(M),$$
i.e. with   initial  distribution $\mu$ the conditional distribution of $X_t$ under $\{\tau>t\}$ remains to be $\mu$ for any $t>0$.
Moreover, we call the Markov process $X_t$ quasi-ergodic with quasi-invariant probability measure  $\mu$, if
$$\lim_{t\to\infty} \E^\nu(f(X_t)|\tau>t)= \mu(f),\ \ \nu\in \scr P_0,$$ where $\nu\in \scr P_0$ means that $\nu\in \scr P$ with  $\E^\nu[1_{\{\tau>t\}}]>0$ for all $t>0$.
In this case, we study the long time behaviors of
$\E^\nu(\W_p(\mu_t,\mu)|\tau>t)$ and $ \W_p(\E^\nu(\mu_t|\tau>t),\mu)$ for $\nu\in \scr P_0$.

\

In the remainder of the paper, we first consider killed diffusion processes on compact manifolds with boundary in Section 2, then study the (reflecting) diffusion processes on compact manifolds in Section 3,
and moreover  investigate general ergodic  Markov processes in Sections 4 and 5 where some typical examples are presented. The corresponding study on weighted empirical measures, time-changed Markov processes, fractional Brownian motions  and non-linear 
Markov processes can be found  in \cite{DU, HT, LW1, LW2, LW3, MT24+,WW, WZ1, Zhu}. 
Finally, we summarize some general results for upper and lower bound estimates on  Wasserstein distance.  

\section{Killed Diffusions on Manifolds}

In this part, we consider the killed diffusion process on a Riemannian manifold with boundary.

Let $M$ be a $d$-dimensional compact connected Riemannian manifold with boundary $\pp M$, let $\d x$ denote the volume measure, and let $N$ be the unit normal vector field of $\pp M$. We call $\pp M$ convex if its second fundamental form is nonnegative,
i.e.
$$\II(U,U):=-\<\nn_U N, U\>\ge 0,\ \ U\in T\pp M.$$
Let $V\in C^1(M)$   such that
$$\mu(\d x):=\e^{V(x)}\d x$$ is a probability measure.
We consider the diffusion process $X_t$ generated by
$$L:=\DD+\nn V$$ which is killed on the boundary $\pp M$, where $\DD$ and $\nn$ stand for the Laplacian and gradient operators respectively.

It is classical that under the Dirichlet boundary condition, the operator  $L$ is a self-adjoint   in $L^2(\mu)$ with discrete spectrum:
$$L\phi_i =-\ll_i \phi_i,\ \ \phi_i|_{\pp M}=0,$$
where $\ll_i\asymp i^{\ff 2 d}$ for $i\ge 1$ and $\{\phi_i\}_{i\ge 0}$ is an orthonormal basis of $L^2(\mu),$ see for instance \cite{Chavel}.

Since  $X_t$ is  the diffusion process on $M$ generated by $L$ with killing boundary condition, its  life time is
$$\tau:=\inf\big\{t\ge 0:\ X_t\in\pp M\big\},$$
and
$$\scr P_0=\big\{\nu\in \scr P:\ \nu(\pp M)<1\big\}.$$
Moreover, the diffusion process is quasi-ergodic with unique quasi-invariant probability measure
 $\mu_0(\d x)=\phi_0^2\d\mu.$
The following result shows that in most cases we have $\W_2(\E^\nu(\mu_t|\tau>t),\mu_0)\asymp t^{-1}$.

\beg{thm}[\cite{W21}]  \label{T1}   Let $\mu_0=\phi_0^2\mu$.  Then for any   $\nu \in \scr P_0,$  
  \beg{align*} &\lim_{t\to\infty} \big\{t^2\W_2(\E^\nu[ \mu_t|\tau>t], \mu_0)^2\big\}\\
&= \ff1 {\{\mu(\phi_0)\nu(\phi_0)\}^2} \sum_{i=1}^\infty \ff{\{\nu(\phi_0)\mu(\phi_i)+ \mu(\phi_0) \nu(\phi_i)\}^2}{(\ll_i-\ll_0)^3}>0, \end{align*}
and the limit is finite if and only if
  $  \nu  \in \D\big((-\hat L)^{-\ff 3 2}\big),$    i.e.
   $$\sum_{i=1}^\infty \nu(\phi_i)^2 \ll_i^{-3}<\infty.$$
 \end{thm}

By the Sobolev embedding theorem, we have $\scr P\subset   \D((-\hat L)^{-\ff 3 2})$ for $d\le 6$, and $\D((-\hat L)^{-\ff 3 2}) \supset L^{\ff{2d}{d+6} }(\mu)$  for $d>6$. Hence, by Theorem \ref{T1},    $\W_2(\E^\nu[ \mu_t|\tau>t], \mu_0)\asymp t^{-1}$  holds if   either  $d\le 6$   or   $d\ge 7$  and
  $\ff{\d\nu}{\d\mu}  \in L^{\ff{2d}{d+6} }(\mu).$

When $d\ge 7$ but $\nu\notin \D((-\hat L)^{-\ff 3 2})$, the limit in the above theorem becomes $\infty$, so that the convergence of $W_2(\E^\nu[ \mu_t|\tau>t], \mu_0)$ is slower than $t^{-1}$, for which the exact convergence rate remains open.

  Next, we   consider  the long time behavior of $\E^\nu(\W_2(\mu_t,\mu_0)|\tau>t).$

\beg{thm}[\cite{WJEMS}]  \label{T2} Let $\nu\in \scr P_0$. The the following assertions hold.
   \beg{enumerate} \item[$(1)$]  When $\pp M$ is convex, we have
   $$ \liminf_{t\to\infty}  \Big\{t\E^\nu\big[\W_2(\mu_{t}, \mu_0)^2\big|t<\tau\big]\Big\}=   \sum_{i=1}^\infty \ff{ 2}{(\ll_i-\ll_0)^2}. $$
   \item[$(2)$]  If  $\pp M$ is non-convex, then there exists a constant $c\in (0,1]$ such that
\beg{align*} &c \sum_{i=1}^\infty \ff{ 2}{(\ll_i-\ll_0)^2}\le   \liminf_{t\to\infty}\Big\{ t\E^\nu\big[\W_2(\mu_{t}, \mu_0)^2\big|t<\tau\big]\Big\}\\
&\le  \limsup_{t\to\infty}\Big\{t \E^\nu\big[\W_2(\mu_{t}, \mu_0)^2\big|t<\tau\big]\Big\}\le   \sum_{i=1}^\infty \ff{ 2}{(\ll_i-\ll_0)^2}. \end{align*}
\end{enumerate}
\end{thm}
 Since $\ll_i-\ll_0\asymp i^{\ff 2 d}$, it is easy to see that
 $$ \sum_{i=1}^\infty \ff{ 2}{(\ll_i-\ll_0)^2}<\infty$$ if and only if $d\le 3$.
 So, Theorem \ref{T2} shows that $\E^\nu\big[\W_2(\mu_{t}, \mu_0)^2\big|t<\tau\big]\asymp t^{-1}$ for $d\le 3$, but the convergence is slower than   than $t^{-1}$ for $d\ge 4$.
 The next result gives the exact convergence rate for $d\ge 4$ where $d=4$ is critical, but for completeness we also include the case for $d\le 3$.
 The lower bound estimate for $d=4$ is due to \cite{W23NS} where more general Markov processes are considered, other estimates  are taken from \cite{WJEMS}.

 \beg{thm}[\cite{WJEMS, W23NS}]   \label{T3}  Let $\nu\in \scr P_0$.  Then  for larger $t>1$
   $$    \E^\nu\big[\W_2(\mu_{t}, \mu_0)^2\big|t<\tau\big]\asymp  \beg{cases} t^{-1}, \ &\text{if}\ d\le 3,\\
t^{-1} \log  t,\ &\text{if}\ d=4,\\
  t^{-\ff 2{d-2}},\ &\text{if}\ d\ge 5. \end{cases} $$
  \end{thm}

\section{(Reflected) Diffusions on Manifolds  }

 Let  $M$ be a  $d$-dimensional connected compact Riemannian manifold without boundary  or with a boundary  $\pp M$.
In this part, we consider the (reflected, if $\pp M$ exists) diffusion process $X_t$ generated by
  $$L:= \DD+ \nn V+ Z,$$
 where  $V\in C^1(M)$  such that
   $\mu(\d x):=\e^{V(x)}\d x$ is a probability measure, and    $Z$  is  a   $C^1$-vector field with
 $${\rm div}_\mu(Z):= Z V+{\rm div}Z=0,$$ where  $Z V:=\<Z,\nn V\>.$
Then $X_t$ is ergodic with unique invariant probability measure $\mu$. We will study the convergence rate of
$\W_p(\mu_t,\mu)$ when $t\to\infty$.

Corresponding   to the killed case where  the Dirichlet eigenvalue  problem is involved,  in the present case we will use the  Neumann eigenvalue problem   if $\pp M$ exists, and  the closed eigenvalue problem otherwise.

Let  $N$ be the inward unit normal vector field on $\pp M$ if the boundary exists.
Consider the Neumann/closed  eigenvalue problem:
  $$L \phi_i=- \ll_i \phi_i, \  \ \mu(\phi_i^2)=1,\  \  \  N \phi_i|_{\pp M}=0,$$
 where   $\phi_0\equiv 1, \ll_0=0$, and the Neumann condition   $N \phi_i|_{\pp M}=0$ applies only when  $\pp M$  exists.
 It is well known that  $\ll_i\asymp i^{\ff 2 d}$ for $i\ge 0$ and $\{\phi_i\}_{i\ge 0} $ is an orthonormal basis of $L^2(\mu),$ see for instance \cite{Chavel}.

 \subsection{Long time behavior of $\E[\W_2(\mu_t,\mu)^2]$}

 Recall that for any  $0\ne f\in L^2(\mu)$ with $\mu(f)=0$,
 $${\bf V}_f:= \int_0^\infty \mu(fP_tf)\d t\in (0,\infty),$$
where  $P_t$  is the   diffusion semigroup of   $X_t$.    By ${\rm div}_\mu(Z)=0,$  $Z\phi_i:=\<Z,\nn\phi_i\>$ satisfies $\mu(Z\phi_i)=0.$ 
The following result is due to \cite{WZ}  for $Z=0$ and \cite{W23NS} otherwise.

 \beg{thm}[\cite{W23NS,WZ}] \label{T4} Let $X_t$ be the (reflected if $\pp M$ exists) diffusion generated by $L=\DD+\nn V+Z$ on $M$. 
 \beg{enumerate} \item[$(1)$]  When $\pp M$ is either empty or convex,
$$\lim_{t\to\infty}   t\E^\nu [\W_2(\mu_t,\mu)^2]=   \sum_{i=1}^\infty\ff 2 {\ll_i^2}\Big(1-\ff 1{\ll_i}{\bf V}_{Z \phi_i}\Big)  $$
holds uniformly in $\nu\in \scr P.$
\item[$(2)$] If $\pp M$ is non-convex,  then  there exists a constant  $c\in (0,1]$ such that
\beg{align*}  &c\sum_{i=1}^\infty\ff 2 {\ll_i^2}\Big(1-\ff 1{\ll_i}{\bf V}(Z \phi_i)\Big)\le  \liminf_{t\to\infty} \inf_{\nu\in \scr P}  t\E^\nu [\W_2(\mu_t,\mu)^2] \\
&\le \limsup_{t\to\infty} \sup_{\nu\in \scr P}   t\E^\nu [\W_2(\mu_t,\mu)^2]\le  \sum_{i=1}^\infty\ff 2 {\ll_i^2}\Big(1-\ff 1{\ll_i}{\bf V}_{Z \phi_i}\Big).\end{align*}
 \end{enumerate}
\end{thm}

 Theorem \ref{T4} shows that a divergence-free perturbation $Z$ accelerates the convergence of $\E[\W_2(\mu_t,\mu)^2]$  with the exact factor  $-\ff 1{\ll_i}{\bf V}_{Z \phi_i}.$
 This fits well to the observation in \cite{GG} that divergence-free perturbations accelerate the convergence in the algorithm of Gibbs measure.

 On the other hand, since $\ll_i\asymp i^{\ff 2 d}$, Theorem \ref{T4} shows that $E^\nu [\W_2(\mu_t,\mu)^2]\asymp t^{-1}$ if and only if $d\le 3$.
 Correspondingly to Theorem \ref{T3}, the following result also present exact convergence rates for $d\ge 4$, where the lower bound estimate for $d=4$ is due to \cite{W23NS} and other estimates are taken from \cite{WZ}.

\beg{thm}[\cite{W23NS, WZ}]  \label{T5}  The following holds for large $t>0$ uniformly in $\nu\in \scr P$:
   $$    \E^\nu\big[\W_2(\mu_{t}, \mu_0)^2 \big]\asymp  \beg{cases} t^{-1}, \ &\text{if}\ d\le 3,\\
t^{-1} \log (t+1),\ &\text{if}\ d=4,\\
  t^{-\ff 2{d-2}},\ &\text{if}\ d\ge 5. \end{cases} $$
\end{thm}

Indeed, for $d=4$ we  have the following renormalization formula.

\beg{thm}[\cite{TWZ}] \label{T6}  When    $d=4$ and $\pp M$  is  empty, there holds 
  $$\lim_{t\to\infty} \sup_{\nu\in \scr P}\bigg|\frac t{\log t}   \E^\nu [\mathbb W_2^2(\mu_t,\mu)] - \frac{{\rm vol}(M)}{8\pi^2} \bigg|=0,$$
where   ${\rm vol}(M)$ is the volume of   $M$. \end{thm}

When $d\ge 5$, it is not clear whether the the following limit  exists or not:
$$\lim_{t\to\infty} t^{\ff 2 {d-2} }\E [\mathbb W_2^2(\mu_t,\mu)]. $$  

\subsection{Long time behavior of $ t\W_2(\mu_{t}, \mu)^2 $ for $d\le 3$}

We first consider the weak convergence of $ t\W_2(\mu_{t}, \mu)^2 $.

\beg{thm}[\cite{WZ}]  \label{T6}
  If  $d\le 3, Z=0$, and   $\pp M$ is either empty or convex, then
  $$ \lim_{t\to\infty}   t\W_2(\mu_{t}, \mu)^2 =  \sum_{i=1}^\infty \ff{2\xi_i^2}{\ll_i^2}\ \text{  in\ law},$$
 where   $\{\xi_i\}$ are i.i.d. standard normal random variables.  \end{thm}

 Next, we consider the convergence of $ t\W_2(\mu_{t}, \mu)^2 $  in $L^q(\P)$ for any $q\ge 1 $   to the following specific process $\Xi(t)$:
  $$\Xi(t):= \sum_{i=1}^\infty \ff{\psi_i(t)}{\ll_i},$$  where
 $$  \psi_i(t):=\ff 1 {\ss t}\int_0^t \phi_i(X_s)\d s.$$
Recall that by the central limit theorem of  \cite{Wu2},
  $\psi_i(t)\to N(0, 2{\bf V}_{\phi_i}) $ weakly with
  $$ {\bf V}_{\phi_i} =\int_0^\infty\mu(\phi_iP_t\phi_i)\d t= \ff 1{\ll_i}\Big(1-\ff 1{\ll_i}{\bf V}_{Z \phi_i}\Big).$$
We will consider initial distributions in the classes
  $$\scr P_{k,R}:=\big\{\nu=\rr_\nu\mu\in \scr P(M):\  \|\rr_\nu\|_{L^k(\mu)}\le R\big\},\ \ k,R\ge 1.$$

 \beg{thm}[\cite{WWZ}] \label{T7} Let  $\pp M$ be empty or convex. If  $d\le 2$, then
   $$\lim_{t\to\infty} \sup_{\nu\in \scr P} \E^\nu \big[\big|t\W_2(\mu_t,\mu)^2- \Xi(t)\big|^q\big]=0,\ \ q\ge 1.$$  If   $d=3$, then
   $$\lim_{t\to\infty} \sup_{\nu\in \scr P_{k,R}} \E^\nu \big[\big|t\W_2(\mu_t,\mu)^2- \Xi(t)\big|^q\big]=0,\ \ k,R, q\ge 1.$$
\end{thm}

\subsection{Moment estimates on $\W_p(\mu_t,\mu)$}
 The following result shows that   the exact convergence rate of $\big(\E[\W_p(\mu_t, \mu)^q]\big)^{\ff 2 q}$ is  uniformly in   $p\ge 1$ and $q>0$.

\begin{thm}[\cite{WWZ}] \label{T8}   Let   $p\in [1,\infty)$ and $ q\in (0,\infty)$.
\beg{enumerate}\item[$(1)$]   If $p\le \ff{2d}{(d-2)^+} \lor \ff{d(d-2)}2,$ where $\ff{2d}{(d-2)^+} =\infty$ when $d\le 2$, then   the following asymptotic formula for large $t>1$ holds   uniformly in $\nu\in \scr P$:
   $$
\big(\E^\nu[\W_p(\mu_t, \mu)^q]\big)^{\ff 2 q} \asymp \begin{cases}
t^{-1}, &\textrm{if } d \leq 3,\\
t^{-1}  \log  t, &\textrm{if } d = 4,\\
t^{-\frac{2}{d - 2}}, &\textrm{if } d \geq 5.\end{cases}$$
\item[$(2)$] If $p> \ff{2d}{(d-2)^+} \lor \ff{d(d-2)}2,$    then for any $k,R\ge 1$ the above asymptotic formula  for large $t>1$  holds uniformly in $\nu\in \scr P_{k,R}.$
\end{enumerate} \end{thm}

\section{Exponential Ergodic Markov processes}

In this part, we consider a general framework of  exponential ergodic Markov process, which will be extended in Section 5 to non-exponentially ergodic setting.

Let $(M,\rr)$ be a length space, which is the basic space for analysis on the Wasserstein space (see \cite{AM,V}), 
 i.e. for any $x,y\in M$, the distance $\rr(x,y)$ can be approximated by the length of curves linking $x$ and $y$.
A typical class of length space is the geodesic space, where the distance $\rr(x,y)$ can be reached by the length of a geodesic curve linking $x$ and $y$.

Let $X_t$ be an ergodic Markov process with unique invariant probability measure $\mu$. We estimate the upper bound of
$\E[\W_p(\mu_t,\mu)^2]$.  Since the empirical measure is usually singular with $\mu$, to apply analysis techniques we need to regularize $\mu_t$ using the following introduced   diffusion process on $M$.

Let $\hat X_t$ be a reversible diffusion process on $M$ with   the same   invariant probability measure $\mu,$  and with $\rr$ as the intrinsic distance. Heuristically,
$\hat X_t$ has symmetric Dirichlet form $(\hat \EE,\D(\hat\EE))$ in $L^2(\mu)$ satisfying
$$ \hat\EE(f,f)= \int_M |\nn f|^2\d\mu,\ \ f\in C_{b,L}(M)\subset \D(\hat\EE),$$
where $C_{b,L}(M)$ is the set of all bounded Lipschitz continuous functions on $M$, and
$$|\nn f(x)|:=\limsup_{y\to x} \ff{|f(y)-f(x)|}{\rr(x,y)},\ \ x\in M.$$
More precisely, we assume that $C_{b,L}(M)$ is a dense subset of $\D(\hat\EE)$ under the   $\hat \EE_1$-norm
$$\|f\|_{\hat\EE_1}:= \ss{\mu(f^2)+\hat\EE(f,f)},$$
and the Dirichlet form restricted on $C_{b,L}(M)$ is formulated as
$$\hat\EE(f,g)=\int_M \GG(f,g)\d\mu,\ \ f,g\in C_{b,L}(M),$$ where
$$\GG: C_{b,L}(M)\times C_{b,L}(M)\to \B_b(M)$$ is a symmetric  local square field (champ de carr\'e), i.e. for
   any $f,g,h\in C_{b,L}(M)$ and $\phi\in C_b^1(\R),$  we have
\beg{align*} &\ss{\GG(f,f)(x)}=|\nn f(x)|:=\limsup_{y\to x} \ff{|f(y)-f(x)|}{\rr(x,y)},\ \ x\in M,\\
&\GG(fg,h)= f\GG(g,h)+ g\GG(f,h),\ \ \ \GG(\phi(f), h)=  \phi'(f) \GG(f,h).\end{align*}
Moreover,  the generator $(\hat L,\D(\hat L))$ satisfies the chain rule
$$\hat L\phi(f)= \phi'(f)\hat Lf +\phi''(f) |\nn f|^2,\ \ \ f\in \D(\hat L)\cap C_{b,L}(M), \phi\in C^2(\R).$$

We make the following assumption.

 \beg{enumerate} \item[$(A)$]   The reference diffusion semigroup    $\hat P_t$  has heat kernel   $\hat p_t$  with respect to   $\mu$, and there exist constants  $\bb,\ll,d,k\in (0,\infty)$ such that
  \beg{align}
      \label{A1} &\|\nn   \hat P_t f\|_{L^2(\mu)\to L^p(\mu)}\le   k\e^{-\ll t} t^{- \bb}, \ t>0,\\
    \label{A3} &\int_M \big(\hat P_t\rr(x,\cdot)^p\big)^{\ff 2 p}(x)\mu(\d x)  \le k t, \ t\in [0,1],\\
  \label{A4} &\int_M  \hat p_t(x,x) \mu(\d x)\le k     t^{-\ff  d 2}, \ t\in (0,1].\end{align}
  Moreover, there exist constants $\theta_1,\theta_2\in (0,\infty)$ such that
  \beg{align} \label{A5}  \|P_t-\mu\|_{L^2(\mu)}  \le \theta_1\e^{-\theta_2 t},\ t\ge 0.\end{align}  \end{enumerate}

Note that conditions  \eqref{A1}-\eqref{A4}  can be verified by choosing a suitable symmetric reference distance $\rr$.   Indeed, for small $\rr$ the Dirichlet form $\hat\EE$ is larger, so that
$\hat P_t$ has better properties.
When  $p=2,$  \eqref{A1}  holds for  $\bb=\ff 1 2$  if and only if 
  $\hat L$ has a spectral gap, i.e.
  $${\rm gap}(\hat L):=\inf\big\{\hat\EE(f,f): \mu(f^2)-\mu(f)^2=1\big\}>0.$$

The only condition on  $X_t$ is the exponential ergodicity  \eqref{A5}, which will be weakened later on by allowing more general ergodic rate.

 Let $K :=    \bb+\ff d 4.   $ 
The convergence rate of   $ \E [\W_p(\mu_t,\mu)^2]$ will be given by
  $$  \xi_K(t):= \beg{cases} t^{-1}, &\text{if}\ K <1,\\
t^{-1}[\log t]^2, &\text{if}\ K =1, \\
t^{-\ff{1}{2K -1 }}, &\text{if}\ K >1. \end{cases}$$

 \beg{thm}[\cite{W24}] \label{T9} Assume $(A)$ and let $K :=    \bb+\ff d 4.   $   Then   there exists a constant  $c>0$ such that for any   $t>0,$
\beq\label{ER}  \E^\mu[\W_p(\mu_t,\mu)^2]\le  c  \xi_K(t),\ \ \ t>0. \end{equation} 
If     $P_t$   has heat kernel    $p_t$  with respect to   $\mu$, then for any $q\in [1,2], t>1$ and  $x\in M,$
  \beg{align*}  \E^x [\W_p(\mu_t,\mu)^q]\le \ &\ff{2^{q-1}}{t^q} \int_0^1 \E^x \big[\mu(\rr(X_s,\cdot)^p)^{\ff q p}\big]\d s
  + 2^{q-1}\|p_1(x,\cdot)\|_{L^{\ff 2{2-q}}(\mu)} (c \xi_K(t-1)^{\ff q 2}.\end{align*}
 In particular, when $p_1$ is bounded, there exists a constant $c>0$ such that 
 \beq\label{ER'}  \sup_{x\in M} \E^x[\W_p(\mu_t,\mu)^2]\le  c  \xi_K(t-1),\ \ \ t>1. \end{equation} 
\end{thm}
Comparing with the exact convergence rate for elliptic diffusions on compact manifolds, the present convergence rate is less sharp. However, as a universal convergence rate for arbitrary exponential ergodic Markov processes,
$\xi_K(t)$ is almost optimal. 
To see this, let us consider the following example.

\paragraph{Example 4.1  (Markov processes on   compact manifolds). }
Let $M$ be a $d$-dimensional compact connected Riemannian manifold possibly with a boundary $\pp M$, let $\mu(\d x)=\e^{V(x)}\d x$ be a probability measure on $M$ for some $V\in C^2(M)$, and let 
 $\hat L:=\DD+\nn V$ (with Neumann boundary condition if $\pp M$ exists). Then  \eqref{A1}-\eqref{A4} hold  for $\bb=\ff 12$ so that $K=   \ff 12+ \ff d 4$.
By Theorem \ref{T9}, 
   for any Markov process on $M$ satisfying \eqref{A5} for some constant $\theta_1,\theta_2>0,$  there exists a constant $c>0$ such that \eqref{ER} holds, and  \eqref{ER'} holds when $p_1$ is bounded, for  
  $$\xi_K(t):= \beg{cases} t^{-1}, &\text{if}\ d <2,\\
t^{-1}[\log t]^2, &\text{if}\ d =2, \\
t^{-\ff{2}{d}}, &\text{if}\ d >2. \end{cases}$$
On the other hand, according to \cite{WW},  for $X_t$ being the $\alpha$-stable time changed process of $\hat X_t$, 
$$\E^x[\W_2(\mu_t,\mu)^2]\asymp \beg{cases} t^{-1}, &\text{if}\ d <2(1+\aa),\\
t^{-1} \log t, &\text{if}\ d =2(1+\aa), \\
t^{-\ff{2}{d-2\aa}}, &\text{if}\ d >2(1+\aa). \end{cases}$$
Since $\xi_K(t)$ is the universal convergence rate for all $\aa\in (0,1)$, it is reached by the exact rate as $\aa\to 0$, except $[\log t]^2$  in the critical case. 

\

Next, we consider a class of Markov processes on $\R^n$.

\paragraph{Example 4.2 (Markov processes on   $\R^n$). }
 Let $M=\R^n$,  let   $V\in C^2(\R^{n})$ such that
$$V(x)=\psi(x) + (1+\theta |x|^2)^\tau,\ \ x\in\R^n, $$ where  $\psi\in C_b^2(\R^{n}),$ $\theta>0,\tau\in (\ff 1 2,\infty]$ are constants.
Let $$\mu(\d x)=\mu_V(\d x):= \ff{\e^{-V(x)}\d x}{\int_{\R^n}\e^{-V(x)}\d x}.$$  Then for  any Markov process on $\R^n$
satisfying \eqref{A5},  there exists a constant $c>0$ such that  
$$ \E^\mu[\W_2(\mu_t,\mu)^2] \le c \beg{cases} t^{-1}, &\text{if}\ n=1, \tau>1,\\
t^{-1} [\log(2+t)]^2, &\text{if}\ n=1, \tau=1,\\
t^{-\ff {2\tau-1}{\tau n}}, &\text{ otherwise}. \end{cases} $$

\beg{proof}   Let $\hat L=\DD-\nn V.$  By Theorem \ref{T9}, it suffices to verify \eqref{A1}-\eqref{A4} for $p=2,\bb=\ff 1 2,$ and $d=\ff{2\tau n}{2\tau-1}.$ Since
$$\lim_{|x|\to\infty} \hat L|\cdot|(x)=-\infty<0,$$
  \cite[Corollary 1.4]{W99} ensures  ${\rm gap}(\hat L)>0$, so that  
    \eqref{A1} holds  for $p=2$ and $\bb=\ff 1 2.$

 Next, by  \cite[Theorem 2.4.4]{Wbook} and $\nn^2 V\ge -c_1I_n,$ where $I_n$ is the $n\times n$-unit matrix, we find a constant $c_2>0$ such that
 \beq\label{AA}\hat p_r(x,x)\le \ff {c_2} {\mu(B(x,\ss r))},\ \ x\in \R^{n}, r\in (0,1],\end{equation}
 where $B(x,r):=\{y\in \R^n: |x-y|<r\}, r>0.$
Then \eqref{A4} with $d=\ff{2\tau n}{2\tau-1}$  follows provided
\beq\label{CC}  \int_{\R^{n }}\ff {\mu(\d x)} {\mu(B(x,r))}\le c r^{-\ff{2\tau n}{2\tau-1}},\ \ r\in (0,1], x\in\R^n\end{equation} holds for some constant $c>0.$
 Since  $\psi$ is bounded, there exists a constant $C>1$ such that
$$C^{-1} \e^{-(1+\theta |x|^2)^\tau}\d x\le \mu(\d x)\le C \e^{-(1+\theta |x|^2)^\tau}\d x.$$ So, \eqref{CC} is trivial for $|x|\le 1$. On the other hand, when $|x|\ge 1$ we have
$$\ff{|x|}2\le |x|-\ff r 4 \le |x|,\ \   r\in [0,1],$$ so   we may find a constant  $c_3>0$ such that
\beg{align*} &\Big(1+ \theta \Big|x- \ff{rx}{4|x|}\Big|^2\Big)^\tau= (1+\theta |x|^2 )^\tau +\int_0^r \ff{\d}{\d s}  \Big(1+ \theta \big(|x|- \ff{s}{4}\big)^2\Big)^\tau\d s\\
&= (1+\theta |x|^2 )^\tau - \ff{\tau\theta}2\int_0^r \Big(1+ \theta \big(|x|- \ff{s}{4}\big)^2\Big)^{\tau-1} \Big(|x|-\ff s 4\Big)\d s\\
&\le  (1+\theta |x|^2 )^\tau- c_3 r |x|^{2\tau-1}.\end{align*}
Hence,  there exist  constants  $c_4, c_5>0$ such that for $|x|\ge 1$ and $r\in (0,1]$,    
\beq\label{SO} \beg{split} &\mu(B(x,r))  \ge c_4 \int_{B\big(x-\ff{r x}{2|x|},\ff r 4\big)}\e^{-(1+\theta  |y|^2)^\tau}\d y\\
&\ge c_5r^n  \e^{-(1+ \theta |x- \ff{rx}{4|x|}|^2)^\tau}  \ge c_5 r^{n}\e^{-(1+\theta |x|^2)^\tau + c_3 r |x|^{2\tau-1}}.\end{split}\end{equation}
 Therefore, there exist  constants $c_6,c_7>0$ such that
\beg{align*}& \int_{\R^{n }}\ff {\mu_V(\d x)} {\mu_V(B(x,r))}\le c_6 r^{-n} \int_{\R^n} \e^{- c_3  r  |x|^{2\tau-1}} \d x \\
& =c_7 r^{-n}\int_0^\infty s^{n-1} \e^{-c_3 rs^{2\tau-1}}\d s= c_7 r^{-\ff{2\tau n}{2\tau-1}}\int_0^\infty s^{n-1} \e^{-c_3 s^{2\tau-1}}\d s.\end{align*}
 Thus, \eqref{CC} holds for some constant $c>0.$

 Finally,   it is easy to see that  $\nn^2 V\ge -c I_n$ and $|\nn V(x)|^2\le c(1+|x|^{4\tau})$ hold for some constant $c>0$. So, we find a constant $c_8>0$ such that
\beq\label{BO} \beg{split} &\hat L|x-\cdot|^2= 2n + 2 \<\nn V, x-\cdot\> = 2n+ 2\<\nn V(x), x-\cdot\> -2\<\nn V(x)-\nn V, x-\cdot\> \\
&\le 2n + |\nn V(x)|^2 + |x-\cdot|^2 + 2c_1|x-\cdot|^2 \le c_8 (1+|x|^{4\tau}+|x-\cdot|^2),\ \ x\in \R^{n}.\end{split} \end{equation} This implies
 \beq\label{B} \hat P_t|x-\cdot|^2(x)=\E^x|x-\hat X_t|^2 \le c_8\big(1+|x|^{4\tau}  |\big)t\e^{c_8t},\ \ x\in\R^{n}, \ t>0.\end{equation}
 Noting that  $\mu( |\cdot|^{4\tau})<\infty,$   we verify condition \eqref{A3} for $p=2$ and some constant $k>0.$
 \end{proof}

As a special case of Example 4.2, we consider the stochastic Hamiltonian system, a typical degenerate SDE  for  $X_t=(X_t^{(1)},X_t^{(2)})$  on   $\R^{n+m}=\R^{n}\times\R^{m}$ ($n,m\ge 1$  may be different):
  \beq\label{E2} \beg{cases} \d X_t^{(1)}= \kk Q X_t^{(2)}\d t,\\
\d X_t^{(2)}= \ss 2 \,\d W_t - \big\{\theta Q^* X_t^{(1)} +\kk X_t^{(2)}\big\}\d t,\end{cases}\end{equation}
where $W_t$ is the $m$-dimensional Brownian motion,  $Q\in\R^{n\otimes m} $ such that $QQ^*$ is invertible, and $\kk,\theta>0$ are constants. 
Let
$$\scr N_\theta(\d x_1):=\Big(\ff{\theta}{2\pi}\Big)^{\ff {n} 2}\e^{-\ff \theta  2 |x_1|^2}\d x_1,\ \ \scr N_\kk(\d x_2):=\Big(\ff{\kk}{2\pi}\Big)^{\ff {m} 2}\e^{-\ff \kk  2 |x_2|^2}\d x_2.$$
By \cite{W17}, where more general degenerate models are considered, the associated Markov semigroup $P_t$ is exponentially ergodic in entropy, hence \eqref{A5} holds. 
So, as shown in Example 4.2, \eqref{A1}-\eqref{A4} hold for $p=2,\bb=\ff 1 2,$ and $d=\ff{2(n+m)}{2\tau-1}.$
 Therefore, for any time-changed process of $X_t$, there exists a constant $c>0$ such that 
 $$ \E^\mu[\W_2(\mu_t,\mu)^2] \le c  
t^{-\ff {1}{m+ n}}.   $$ 
 
\section{More general  ergodic Markov processes}

For some infinite-dimensional models, see for instance \cite{WAAP}, \eqref{A4} fails for any $d\in (0,\infty)$, but there may be  a decreasing function $\gg: (0,\infty)\to (0,\infty)$ such that
$$ \int_M \hat p_t(x,x)\mu(\d x)\le \gg(t),\ \ t>0.$$
Next, in infinite-dimensional case the condition \eqref{A3} may be invalid for small time, see Corollary \ref{5.2} below. 
Moreover, in case that $P_t$ is not $L^2$-exponential ergodic, by the weak Poincar\'e inequality which holds for a broad class of ergodic Markov processes,   see \cite{RW01}, we  have 
$$\lim_{t\to\infty}  \|P_t-\mu\|_{L^\infty(\mu)\to L^2(\mu)}= 0.$$
To cover these  situations for which Theorem \ref{T9} does not apply, we present the following result for the empirical measure $\mu_t$ of the Markov process $X_t$ with semigroup $P_t$.

\beg{thm}[\cite{W24}]\label{TE} Assume   $\eqref{A1}$,  $\eqref{A3}$. If there exist a constant $q\in [1,\infty],q'\in[\ff q {q-1},\infty],$  a decreasing function $\gg: (0,\infty)\to (0,\infty)$ and an increasing continuous function $h: [0,1]\to [0,\infty)$ with $h(0=0$ such that   
\beg{align} \label{A3'} & \int_M \big(\hat P_t\rr(x,\cdot)^p\big)^{\ff 2 p}(x)\mu(\d x)  \le h(t),\ \ t\in (0,1], x\in M, \\
 \label{A5'} &\lim_{t\to \infty}  \|P_t-\mu\|_{L^{q'}(\mu)\to L^{\ff{q}{q-1}}(\mu)} =0,\\
 \label{A4'} & \int_M \|\hat p_{\ff r 2}(y,\cdot)\|_{L^q(\mu)}\|\hat p_{\ff r 2}(y,\cdot)\|_{L^{q'}(\mu)} \mu(\d y)\le \gg(r),\ \ \ r>0.\end{align}
 For any $t>0$, let
 $$\xi(t):=  \inf_{r\in (0,1]} \bigg\{\ff {\int_0^t \|P_{s}-\mu\|_{L^{q'}(\mu)\to L^{\ff q {q-1}}(\mu)}\d s} t \bigg(\int_0^1 \ff{\ss{\gg(r+s)}}{(r+s)^\bb}\d s\bigg)^2 +h(r)\bigg\}.$$
Then there exists a constant $c>0$ such that for any $t>0,$
\beq\label{TE1} \E^\mu[\W_p(\mu_t,\mu)^2] \le c \xi(t).\end{equation}
If    $P_t$   has heat kernel    $p_t$  with respect to   $\mu$, then for any $q\in [1,2], t>1$ and  $x\in M,$
  \beq\label{TE2}   \E^x [\W_p(\mu_t,\mu)^q]\le c t^{-q} \int_0^1 \E^x \big[\mu(\rr(X_s,\cdot)^p)^{\ff q p}\big]\d s
  + c\big( \|p_1(x,\cdot)\|_{L^{\ff 2{2-q}}(\mu)}   \xi(t-1)\big)^{\ff q 2}.\end{equation}
 \end{thm}  
 To verify Theorem \ref{TE}, we present below a simple example where $P_t$ only has algebraic convergence in $\|\cdot\|_{L^\infty(\mu)\to L^2(\mu)},$ so Theorem \ref{T9} does not apply. 

\paragraph{Example 5.1.}  Let $M=[0,1], \rr(x,y)=|x-y|$ and $\mu(\d x)= \d x.$ For any $l\in (2,\infty),$ let $X_t$ be the diffusion process on $M\setminus \{0,1\}$ generated by
$$L:= \big\{x(1-x)\big\}^l\ff{\d^2}{\d x^2} + l \big\{x(1-x)\big\}^{l-1}(1-2x) \ff{\d}{\d x}.$$
Then  Theorem \ref{T9} does not apply, but by Theorem \ref{TE} there exists a constant $c>0$ such that for any $t>0$, 
\beq\label{CV} \E^\mu[\W_p(\mu_t,\mu)^2]\le c \beg{cases} t^{-1}, &\text{if} \ l\in (2,5), p \in [2, \ff{13-l}4), \\
t^{-1}[\log(2+t)]^3,  &\text{if} \ l\in (2,5], p=\ff{13-l}4, \\
 \big[t^{-1} \log(2+t)]^{\ff 8{4p+l-5}}, &\text{if} \   l\in (2,5], p>\ff{13-l}4,\\
 t^{-\ff 4{l-1}}[\log(2+t)]^2, &\text{if}\ l>5, p=2\\
 t^{-\ff 8 {p(l-1)}},  &\text{if} \ l>5, p>2.\end{cases}\end{equation} 

\beg{proof}  We first observe that \eqref{A5} fails, so that Theorem \ref{T9} does not apply. Indeed, the Dirichlet form of $L$ satisfies
\beq\label{DF} \EE(f,g)=\int_0^1 \big\{x(1-x)\big\}^l (f'g')(x)\d x,\ \ f,g\in C_b^1(M)\subset \D(\EE).\end{equation} 
Let $\rr_L$ be the  intrinsic distance  function to the point $\ff 1 2\in M.$ We find a constant $c_1>0$ such that
$$\rr_L(x)= \bigg|\int_{\ff 1 2}^x \big\{s(1-s)\big\}^{-\ff l 2}\d s\bigg|\ge c_1 \big(x^{1-\ff l 2}+ (1-x)^{1-\ff l 2}\big),\ \ x\in M.$$
Then  for any $\vv>0,$ we have $\mu(\e^{\vv \rr_L})=\infty,$ so that by \cite{AMS},  ${\rm \gap}(L)=0.$ 
On the other hand, since $L$ is symmetric in $L^2(\mu)$, by \cite[Lemma 2.2]{RW01},    \eqref{A5} implies  the same inequality for $k=1$, so that ${\rm gap}(L)\ge \ll>0$.  Hence,  \eqref{A5} fails. 

To apply Theorem \ref{TE}, let $\hat P_t$ be the standard Neumann heat semigroup on $M$ generated by $\DD$. 
It is classical that \eqref{A1} and  \eqref{A3'}   with $h(t)=kt$  hold  for some constant $k>0$ and 
\beq\label{BMM} \bb= \ff 1 2+\ff {p-2}4.\end{equation} Moreover,  there exists a constant $c_2>1$ such that 
$$\|\hat P_{\ff r 2}\|_{L^m(\mu)\to L^n(\mu)}\le c_2(1+r^{-\ff {n-m}{2nm}}),\ \ \ 1\le m\le n\le\infty, r>0,$$ so that for $q'= \infty$ and $q>1$, 
\beg{align*}  &\|\hat p_{\ff r 2}(y,\cdot)\|_{L^q(\mu)}\|\hat p_{\ff r 2}(y,\cdot)\|_{L^{q'}(\mu)} \mu(\d y)\le c_2 \|\hat P_{\ff r 2}\|_{L^1(\mu)\to L^q(\mu)}  (1+ r^{-\ff 1 2}) \\
&\le c_2^2 (1+r^{-\ff {q-1} {2q}}) (1+r^{-\ff 1 2}),\ \ r>0.\end{align*}
Hence, there exists a constant $c_3>0$ such that \eqref{A4'} holds for
$$\gg(r)= c_3  (1+ r^{-\ff{2q-1}{2q}}).$$ Combining this with \eqref{BMM}, we find a constant $k>0$ such that for any $r\in (0,1),$ 
\beq\label{GMM} \eta(r):= \bigg(\int_0^1 \ff{\ss{\gg(r+s)}}{(r+s)^\bb}\d s\bigg)^2\le k\beg{cases}  
1, &\text{if}\ 1<q<\ff 1 {p-2},\\
[\log(1+r^{-1})]^2, &\text{if}\ 1<q=\ff 1 {p-2},\\
r^{\ff{1-(p-2)q}{2q}}, &\text{if}\ q>  1\lor\ff{1}{p-2}.\end{cases}\end{equation}

To calculate $\|P_t-\mu\|_{L^{q'}(\mu)\to L^{\ff q{q-1}}(\mu)} $  for $q'=\infty$, we apply  the weak Poincar\'e inequality studied in \cite{RW01}. Let
$$M_s=[s,1-s],\ \ \ s\in (0,1/2).$$
Noting that $\mu(\d x)=\d x$ and letting  $\nu(\d x)= \{x(1-x)\}^l]\d x,$ we find a constant $c_4>0$ such that 
$$\sup_{r\in [s,1/2]} \mu([r,1/2])  \nu([s,r]) \le  2^l \sup_{r\in [s,\ff 1 2]} \Big(\ff 1 2 -r\Big)  \big(s^{1-l}- r^{1-l} \big) 
\le c_4 s^{1-l},\ \ s\in (0,1/2).$$ By the weighted Hardy inequality,  see for instance \cite[Proposition 1.4.1]{Wbook},  
we have
$$\mu(f^2 1_{[s,\ff 1 2]})\le 4 c_4 s^{1-l}\nu(|f'|^2),\ \ f\in C^1([s,1/2]), f(1/2)=0.$$
By symmetry, the same holds for $[\ff 1 2,1-s]$ replacing $[s,\ff 1 2].$ So, according to \cite[Lemma 1.4.3]{Wbook},   we derive
$$\mu(f^2 1_{M_s}) \le 4 c_4 s^{1-l} \nu(|f'|^21_{M_s}) +\mu(f1_{M_s})^2,\ \ f\in C^1([s,s-1]).$$
Combining this with \eqref{DF}, for any $f\in C^1_b(M)$ with $\mu(f)=0$, we have $\mu(f1_{M_s})= -\mu(f1_{M_s^c})$ so that 
\beg{align*} &\mu(f^2) = \mu(f^21_{M_s^c}) + \mu(f^21_{M_s})\le \mu(f^21_{M_s^c}) +  4c_4 s^{1-l} \EE(f,f)  +\mu(f1_{M_s^c})^2,\\
&\le 4 c_4 s^{1-l}   \EE(f,f) + 2 \|f\|_\infty^2 \mu(M_s^c)^2\le 4 c_4 s^{1-l}\EE(f,f) + 8 s^2 \|f\|_\infty^2,\ \ s\in (0, 1/2).\end{align*}
For any $r\in (0,1),$ let $s= (r/8)^{\ff 1 2}.$ We find a constant $c_5>0$ such that 
$$\mu(f^2) \le c_5 r^{-\ff {l-1} 2} \EE(f,f) +r\|f\|_\infty^2,\ \ r\in (0,1), \mu(f)=0, f\in C^1_b(M).$$ 
By \cite[Corollary 2.4(2)]{RW01}, this implies 
$$\|P_t-\mu\|_{L^\infty(\mu)\to L^2(\mu)}=\|P_t-\mu\|_{L^2(\mu)\to L^1(\mu)}\le c_5 (1+t)^{-\ff 2 {l-1}},\ \ t>0$$ for some constant $c_5>0.$ 
Since $ P_t $ is contractive in $L^n(\mu)$ for any $n\ge 1$,  this together with the interpolation theorem  implies 
$$\|P_t-\mu\|_{L^\infty(\mu)\to L^{\ff q {q-1}}(\mu)}\le c_6 (1+t)^{-\ff {4(q-1)}{q(l-1)}},\ \ t>0.$$
Noting that $q'=\infty,$ we find a constant $k>0$ such that
\beq\label{GGT} \GG(t):= \ff 1 t \int_0^t \|P_s-\mu\|_{L^{q'}(\mu)\to L^{\ff q {q-1}}(\mu)}\d s\le k\beg{cases} 
 t^{-1},  &\text{if}\  l\in (2,5), q>\ff 4{5-l},   \\
t^{-1} \log(2+t),  &\text{if}\ l=5, q=\infty,\\
(1+t)^{-\ff{4 }{l-1}},  &\text{if}\ l>5, q=\infty.\end{cases}\end{equation} 
We now prove the desired estimates case by case.

(1) Let $l\in (2,5)$ and $p\in [2,\ff{13-l}4).$ Taking  $q\in (\ff 4{5-l}, \ff 1 {p-2})$ in \eqref{GMM} and \eqref{GGT},  we obtain 
$$\inf_{r\in (0,1]} \big\{\eta(r) \GG(t)+r\big\}\le k \inf_{r\in (0,1]}\big\{t^{-1}+r\big\}= kt^{-1}.$$
So, the desired estimate follows from Theorem \ref{TE}.

(2) Let $l\in (2,5]$ and $p=\ff{13-l}4.$ Taking $q= \ff 4{5-l}=\ff 1 {p-2}$ in  \eqref{GMM} and \eqref{GGT} we find a constant $c>0$ such that 
$$\inf_{r\in (0,1]} \big\{\eta(r) \GG(t)+r\big\}\le k \inf_{r\in (0,1]}\big\{t^{-1}[\log(2+t)][\log (1+r^{-1})]^2+r\big\}\le c t^{-1}[\log (2+t)]^3.$$
This implies the desired estimate  according to  Theorem \ref{TE}.

(3) Let $l\in (2,5]$ and $p>\ff{13-l}4.$ We have  $q:= \ff 4 {5-1}>\ff 1 {p-2},$ so that   \eqref{GMM} and \eqref{GGT} imply
$$\inf_{r\in (0,1]} \big\{\eta(r) \GG(t)+r\big\}\le k \inf_{r\in (0,1]}\big\{t^{-1}[\log(2+t)] r^{-\ff{4p+l-13} 8}+r\big\}\le c \big[t^{-1}\log (2+t)\big]^{-\ff{4p+l-5}8} $$
for some constant $c>0$, which implies the desired estimate by Theorem \ref{TE}. 

(4) Let $l>5$ and $p=2$. By taking $q=\infty$ in  \eqref{GMM} and \eqref{GGT}, we find a constant $c>0$ such that  
$$\inf_{r\in (0,1]} \big\{\eta(r) \GG(t)+r\big\}\le k \inf_{r\in (0,1]}\big\{t^{-\ff 4{l-1}}[ \log(1+r^{-1})]^2 +r\big\}\le c  t^{-\ff 4{l-1}}[\log(2+t)]^2.$$
By Theorem \ref{TE},  the desired estimate holds. 

(5) Let $l>5$ and $p>2$. By taking $q=\infty$ we find a constant $c>0$ such that  \eqref{GMM} and \eqref{GGT} imply 
$$\inf_{r\in (0,1]} \big\{\eta(r) \GG(t)+r\big\}\le k \inf_{r\in (0,1]}\big\{t^{-\ff 4{l-1}} r^{-\ff{p-2}2} +r\big\}\le c  t^{-\ff 8{p(l-1)}}$$
for some constant $c>0.$ Hence 
   the desired estimate holds according to Theorem \ref{TE}. 
 
\end{proof} 

Finally, we consider  the following semilinear SPDE  on a separable Hilbert space $\H$: Consider the following SDE on a separable Hilbert space $\H$:
\beq\label{E1} \d \hat X_t= \big\{\nn V(\hat X_t)-A \hat X_t \big\}\d t + \ss 2\, \d W_t,\end{equation}
where $W_t$ is the cylindrical Brownian motion on $\H$, i.e. 
$$W_t= \sum_{i=1}^\infty B_t^i e_i,\ \ t\ge 0$$ for an orthonormal basis $\{e_i\}_{i\ge 1}$ of $\H$ and a sequence of independent one-dimensional Brownian motions $\{B_t^i\}_{i\ge 1},$ 
 $(A,\D(A))$ is a positive definite self-adjoint operator  and $V\in C^1(\H)$ satisfying the following assumption.
 
\beg{enumerate} \item[$(H_1)$] $A$ has discrete spectrum with eigenvalues $\{\ll_i>0\}_{i\ge 1}$ listed in the increasing order counting multiplicities  satisfying
$\sum_{i=1}^d \ll_i^{-\dd}<\infty$ for some constant $\dd\in (0,1)$, and $V\in C^1(\H)$, $\nn V$ is Lipschitz continuous in $\H$ such that \beq\label{VV} \<\nn V (x)- \nn V(y), x-y\>\le (K+\ll_1)|x-y|^2,\ \ x,y\in \H\end{equation}
holds for some constant $K\in\R$. Moreover,
$Z_V:= \mu_0(\e^V)<\infty,$  where
  $\mu_0$ is the centered Gaussian measure on $\H$ with covariance operator $A^{-1}$. 
\item[$(H_2)$]  There exists  an increasing function $\psi: (0,\infty)\to [0,\infty)$ such that
$$ |V(x)|\le \ff 1 2 \big(\psi(\vv^{-1})+  \vv   |x|^2\big),\ \ x\in \H,\vv>0.$$
\item[$(H_3)$] There exist   constants $c>0$ and   $\theta\in [0, \ll_1)$
$$ \ |\nn V(x)|\le c+ \theta |x|,\ \ x\in \H.$$  \end{enumerate}

Under $(H_1)$,  for any $\F_0$-measurable random variable $X_0$ on $\H$, \eqref{E1} has a unique mild solution, and there exists an increasing function $\psi: [0,\infty)\to (0,\infty)$ such that
$$ \E[\|\hat X_t\|_\H^2]\le  \psi(t)\big(1+ \E[|X_0|^2]\big),\ \ t\ge 0,$$
see for instance \cite[Theorem 3.1.1]{W13}.

 Let  $\hat P_t$ be the associated Markov semigroup, i.e.
$$\hat P_tf(x):= \E^x [f(\hat X_t)],\ \ t\ge 0, f\in \B_b(\H),\ \ x\in \H,$$
where $\B_b(\H)$ is the class of all bounded measurable functions on $\H$, and $\E^x$ is   the expectation   for the solution $X_t$ of \eqref{E1} with $X_0=x$. In general, for a probability measure $\nu$ on $\H$, let $\E^\nu$ be the expectation for $X_t$ with initial distribution $\nu$.

By $(H_1),$ we define the probability measure
  $$\mu(\d x):= Z_V^{-1} \e^{V(x)}\mu_0(\d x).$$
Then  $\hat P_t$   is symmetric in $L^2(\mu)$.  

\beg{cor} \label{5.2} Assume $(H_1)$ and $(H_2)$.  Let
$$\xi(t):=\inf_{r\in (0,1)}    \bigg( \ff 1 t r^{-1}\e^{kr^{-1}+\psi(kr^{-1})}  +r^{1-\dd} \bigg),\ \ t>0.$$ Let  $X_t=\hat X_{S_t}$ for an increasing stable process $S_t$ with $S_0=0$ which is independent of $X$. 
  \beg{enumerate}
\item[$(1)$] There exists a constant $c >0$ such that
\beq\label{B1} \E^\mu\big[\W_2(\mu_t,\mu)^2\big]\le c\xi(t),\ \ t>0.\end{equation}
\item[$(2)$] If  $(H_3)$ holds,  then there exists a constant $k>0$ such that
\beq\label{B2} \big(\E^x[\W_2(\mu_t,\mu)]\big)^2\le k  \e^{k|x|^2} \xi(t-1),\ \ t> 1, x\in\H. \end{equation}\end{enumerate}
\end{cor}
 
 \beg{proof} By (2.22) in \cite{WAAP},    $(H_1)$ and $(H_2)$   imply \eqref{A4'} for
 $$\gg(r):= \e^{\psi(k\vv^{-1})+ k\vv^{-1}}$$ for some constant $k>0$. As explained after (2.4) in \cite{WAAP} that this implies \eqref{A1} for $p=2$ and $\bb=\ff 1 2$, as well as   
 $$\|P_t-\mu\|_{L^2(\mu)}\le \e^{-\ll_0 t},\ \ t>0$$ for some constant $\ll_0>0$. Moreover, since $k:= \sum_{i=1}^\infty \ll_i^{-\dd}<\infty$, 
 $$\sum_{i=1}^\infty \int_0^t \e^{-2\ll_i(t-s)}\d s \le \sum_{i=1}^\infty \ll_i^{-\dd}t^{1-\dd}=kt^{1-\dd},$$ 
 so that  \eqref{A3'}  holds for $p=2$ and $h(t)=kt^{1-\dd}$. Therefore, \eqref{B1} follows form  \eqref{TE1}. 
 If moreover $(H_3)$ holds, then as shown in the proof of \cite[Corollary 2.2(2)]{WAAP}, we have $p_1(x,x)\le c\e^{c|x|^2}$ for some constant $c>0$, so that 
 \eqref{B2} is implied by \eqref{TE2}. 
 \end{proof} 

The following example extends  \cite[Example 2.1]{WAAP}. 

\paragraph{Example 5.2.} Assume $(H_1) $ and that such that $ V$ is Lipschitz continuous, and let  $X_t=\hat X_{S_t}$ for an increasing stable process $S_t$ with $S_0=0$ which is independent of $X$.    Then $(H_2)$ holds for 
$ \psi(s)=c_0 s$ for some constant $c_0>0$. So, by taking $r= N(\log t)^{-1}$ for a large enough constant $N>0$, we find a constant $c>0$ such that  \eqref{B1} and \eqref{B2}  imply that for large $t>0$ 
\beg{align*}& \E^\mu\big[\W_2(\mu_t,\mu)^2\big]\le c_1 (\log t)^{\dd-1},\\
&  \big(\E^x [ \W_2(\mu_t,\mu)]\big)^2 \le c\e^{c|x|^2} 
(\log t)^{\dd-1}, \ \ x\in\H.  \end{align*}

\section{Some General Estimates on Wasserstein Distance} 
In this part, we introduce some useful techniques in estimating the Wasserstein distance for empirical measures of diffusion processes. 

Let $(M,\rr)$ be a length space, and recall that for a Lipschitz continuous function $f$
$$|\nn f|(x):=\limsup_{y\to x} \ff{|f(x)-f(y)|}{\rr(x,y)}.$$
Let $(\hat L,\D(\hat L))$ be the self-adjoint Dirichlet operator in $L^2(\mu)$ with Dirichlet form $(\hat\EE,\D(\hat\EE)$ satisfying  $C_{b,L}(M)\subset \D(\hat\EE)$ and      
  $$\hat\EE(f,f)=\mu(|\nn f|^2),\ \ \ f\in C_{b,L}(M).$$
  Let $X_t$ be the diffusion process generated by $$L:=\hat L+Z,$$ where 
  $Z$ is a bounded vector field with ${\rm div}_\mu(Z)=0$, i.e.
 \beg{align*} &  \|Z\|_\infty:=\sup_{\|f\|_{Lip}\le 1} |Zf|<\infty,\\
 &\int_M Zf \d\mu=0,\ \ \ f\in C_{b,L}(M).\end{align*}

  \subsection{Upper bound estimate}
   
 According to \cite[Theorem 2]{Ledoux}, for any probability density $f$ of $\mu$, we have
 \beq\label{APP1}  \W_p(f\mu, \mu)^p \le p^p  \mu\big( |\nn(-\hat L)^{-1} (f-1)|^p \big),\ \ p\in [1,\infty).\end{equation}  The idea of the proof goes back to \cite{AMB}, in which   the following estimate is presented
 for probability density functions $f_1,f_2$:
 \beq\label{APP2} \W_2(f_1\mu_1,f_2\mu_2)^2\le \int_M \ff{|\nn (-\hat L)^{-1} (f_2-f_1)|^2}{\scr M(f_1,f_2)}\d\mu,\end{equation}
 where $\scr M(a,b):=1_{\{a\land b>0\}} \ff{\log a-\log b}{a-b}$ for $a\ne b$, and $\scr M(a,a)=1_{\{a>0\}}a^{-1}$.  In general,    for $p\ge 1$,   denote $\scr M_p=\scr M$ if $p=2$, and when $p\ne 2$ let
 $$\scr M_p(a,b)=1_{\{a\land b>0\}} \ff{a^{2-p}-b^{2-p}}{(2-p)(a-b)}\text{\ for\ }a\ne b,\ \ \scr M_p(a,a)= 1_{\{a>0\}}a^{1-p}.$$
 The following result extends  \eqref{APP1} and \eqref{APP2}.

\beg{thm} \label{TA1}  For any   probability  density functions $ f_1$ and $f_2$ with respect to $\mu$ such that  $f_1\lor f_2>0$,
\beg{align*}  \W_p(f_1\mu, f_2\mu)^p \le \min\bigg\{ &p^p 2^{p-1} \int_M \ff{| \nn (-\hat L)^{-1} (f_2-f_1)|^p}{(f_1+f_2)^{p-1}}\d\mu,\  p^p\int_M \ff{| \nn (-\hat L)^{-1} (f_2-f_1)|^p}{f_1^{p-1}}\d\mu,\\
&\qquad \int_M \ff{| \nn (-\hat L)^{-1}(f_2-f_1)|^p }{\scr M_p(f_1,f_2)}\d\mu\bigg\}.\end{align*}

 \end{thm}

\beg{proof}   It suffices to prove for $p>1$.
   Consider the Hamilton-Jacobi semigroup $(Q_t)_{t>0}$ on $C_{b,L}(M)$:
 $$ Q_t \phi:= \inf_{x\in M} \Big\{\phi(x)+ \ff 1 {p t^{p-1}} \rr(x,\cdot)^p\Big\},\ \ t>0, \phi\in C_{b,L}(M).$$
 Then for any $\phi\in C_{b,L} (M)$, $Q_0\phi:= \lim_{t\downarrow 0} Q_t\phi=\phi$, $\|\nn Q_t\phi\|_\infty$ is locally bounded in $t\ge 0$, and $Q_t \phi$ solves the  Hamilton-Jacobi equation
 \beq\label{HK0} \ff{\d}{\d t} Q_t \phi= -\ff {p-1}p  |\nn Q_t\phi|^{\ff p{p-1}},\ \ t>0.\end{equation}
Let $q=\ff{p}{p-1}.$ For any $f\in C_b^1(M)$,  and any increasing function $\theta\in C^1((0,1))$ such that $\theta_0:=\lim_{s\to 0}\theta_s=0, \theta_1:=\lim_{s\to 1}\theta_s=1$,
by   \eqref{HK0} and the integration by parts formula, we obtain
\beg{align*} &\mu_1(Q_1f)- \mu_2(f)= \int_0^1\Big\{\ff{\d}{\d s}  \mu\big([f_1+\theta_s(f_2-f_1)]Q_s f\big)\Big\}\d s \\
& = \int_0^1 \d s \int_M \Big\{\theta_s'(f_2-f_1) Q_s f -\ff{f_1+\theta_s(f_2-f_1)}q |\nn Q_s f|^q\Big\}\d\mu \\
&= \int_0^1  \Big[ \theta_s'\hat\EE \big(  (-\hat L)^{-1}(f_2-f_1),    Q_sf\big) -\mu\Big(\ff{f_1+\theta_s(f_2-f_1)}q |\nn Q_s f|^q\Big)\Big]\d s \\
&\le \int_0^1 \mu \Big(\theta_s'\big|\nn  (-\hat L)^{-1}(f_2-f_1)\big|\cdot\big|\nn   Q_sf\big| - \ff{f_1+\theta_s(f_2-f_1)}q |\nn Q_s f|^q\Big) \d s.  \end{align*}
Combining this with    Young's inequality $ab\le a^p/p+ b^q/q$ for $a,b\ge 0$, we arrive at
\beq\label{SPP} \mu_1(Q_1f)- \mu_2(f)\le \ff 1 p \int_M|\nn (-\hat L)^{-1}(f_2-f_1)|^p \d\mu  \int_0^1 \ff{|\theta_s'|^p}{[f_1+\theta_s(f_2-f_1)]^{p-1}}\d s.\end{equation} 
 By Kantorovich duality formula
$$\ff 1 p\W_p(\mu_1,\mu_2)^p= \sup_{f\in C_b^1(M)} \big\{\mu_1(Q_1 f)-\mu_2(f)\big\},$$
and noting that
\beg{align*} &f_1+\theta_s(f_2-f_1)= f_1+ f_2 - \theta_sf_1 -(1-\theta_s) f_2 \\
&= (f_1+ f_2)\Big(1-\ff{\theta_sf_1}{f_1+ f_2}- \ff{(1-\theta_s) f_2}{f_1+f_2}\Big)\\
&\ge (f_1+ f_2) \min\{1-\theta_s, \theta_s\},\end{align*}
we deduce from \eqref{SPP} that 
\beq\label{ECC} \W_p(\mu_1,\mu_2)^p \le  \int_0^1\ff{|\theta_s'|^p}{\min\{\theta_s, 1-\theta_s\}^{p-1}}\d s \int_M \ff{|\nn (-\hat L)^{-1}(f_1-f_2)|^p}{(f_1+f_2)^{p-1}}\d\mu.\end{equation}
By taking
$$\theta_s= 1_{[0,\ff 1 2]}(s) 2^{p-1} s^p + 1_{(\ff 1 2,1]}(s) \big\{1- 2^{p-1} (1-s)^p\big\}, $$
which satisfies
$$\theta_s'= p 2^{p-1}\min\{s,1-s\}^{p-1},\ \ \min\{\theta_s, 1-\theta_s\}= 2^{p-1}\min\{s,1-s\}^{p},$$ we deduce from \eqref{ECC}
that
 $$\W_p(f_1\mu, f_2\mu)^p \le p^p  2^{p-1} \int_M \ff{| (-\hat L)^{-\ff 1 2 } (f_2-f_1)|^p}{(f_1+f_2)^{p-1}}\d\mu. $$
Next, \eqref{ECC} with $\theta_s= 1-(1-s)^p$ implies
$$\W_p(f_1\mu, f_2\mu)^p \le p^p \int_M \ff{| (-\hat L)^{-\ff 1 2 } (f_2-f_1)|^p}{f_1^{p-1}}\d\mu.$$
Finally, with $\theta_s=s$ we deduce from \eqref{SPP}  that
$$\W_p(f_1\mu, f_2\mu)^p \le \int_M \ff{|(-\hat L)^{-\ff 1 2}(f_2-f_1)|^p }{ \scr M_p(f_1,f_2)}\d\mu.$$
Then the proof is finished.
\end{proof}

We now apply Theorem \ref{TA1} to the   regularized   empirical measure  
$$\mu_{t,\vv}:=  \hat P_\vv^*\mu_t$$ for suitable choice of $\vv=\vv_t\downarrow  0$ as $t\uparrow\infty$.
To  this end,  we make the following assumption on the reference diffusion process $\hat X_t$ introduced in Section 4. 

\beg{enumerate}\item[{\bf (A)}]   The following conditions hold for some $ d\in [1,\infty)$ and  an increasing function $K: [2,\infty)\to (0,\infty)$.
\item[$\bullet$] {\bf  Nash inequality.} There exists a constant $C>0$ such that
\beq\label{N}  \mu(f^2)\le C   \hat \EE (f,f)^{\ff d{d+2}}\mu(|f|)^{\ff 4{d+2}},\ \ f\in\D_0:=\big\{f\in\D( \hat \EE):\ \mu(f)=0\big\}. \end{equation}
\item[$\bullet$]  {\bf Continuity of symmetric diffusion.} For any $p\in [2,\infty)$,
\beq\label{CT}   \E^\mu [\rr(\hat X_0, \hat X_t)^p] =   \int_{  M\times M}  \rr(x,y)^p \hat p_t(x,y) \mu(\d x) \mu(\d y)\le K(p) t^{\ff p 2},\ \ t\in [0,1],\end{equation}
where $\hat p_t$ is the heat kernel of $\hat P_t$ with respect to $\mu$.
\item[$\bullet$] {\bf  Boundedness of Riesz transform.} For any $p\in [2,\infty)$,
\beq\label{RZ} \|\nn (-\hat L)^{-\ff 1 2} f\|_{L^p(\mu)}\le K(p)\|f\|_{L^p(\mu)},\ \ f\in L^p(\mu) \text{ with } \mu(f)=0.\end{equation}
 \end{enumerate}

It is well known that {\bf (A)} holds for the  (reflecting) diffusion process generated by $\hat L:=\DD+\nn V$ considered in Introduction.

Besides the elliptic diffusion process on compact manifolds,  some  criteria on the Nash inequality \eqref{N} are available in  \cite[Section 3.4]{Wbook}. In general,  \eqref{N} implies that for some constant $c_0 >0$,
 \beq\label{HS} \|\hat P_t-\mu\|_{L^p(\mu)\to L^q(\mu)}\le c_0 (1\wedge t)^{-\ff {d(q-p)} {2pq}}\e^{-\ll_1 t},\ \ t>0, \ 1 \leq p \leq q \leq \infty,\end{equation}
 and that $-\hat L$ has purely discrete spectrum with all eigenvalues $\{\ll_i\}_{i\ge 0}$, which are listed in the increasing order counting multiplicities, satisfy
  \beq\label{EG} \ll_i\ge c_1 i^{\ff 2 d},\ \ i\ge 0,\end{equation} for some constant $c_1>0$. The Markov semigroup $\hat P_t$ generated by $\hat L$ has symmetric heat kernel $\hat p_t$ with respect to $\mu$ formulated as
\beq\label{HS2} \hat p_t(x,y)= 1+ \sum_{i= 1}^\infty \e^{-\ll_i t} \phi_i(x)\phi_i(y),\ \ t>0, \, x,y\in M. \end{equation}
All these assertions can be found  for instance in \cite{Davies}.

The condition \eqref{CT} is natural for diffusion processes due to the growth property $\E|B_t-B_0|^p \le c t^{\ff p2}$ for  the Brownian motion $B_t$. There are plentiful results on the boundedness  condition \eqref{RZ} for the  Riesz transform, see \cite{B87, CD11, CTW23} and references therein.

The following result shows that  under assumption {\bf (A)}, the convergence rate of $(\E \W_p^q(\mu_T,\mu))^{\ff 1 q}$ is  given by
 $$\gg_d(t):= \beg{cases} t^{-\ff 1 2 }, &\textrm{if } \ d\in [1,4),\\
t^{-\ff 1 2}  \ss{\log t}, &\textrm{if } \ d = 4,\\
t^{-\frac{1}{d - 2}}, &\textrm{if }\  d\in (4,\infty).\end{cases} $$

\beg{thm}[\cite{WWZ}] Assume {\bf (A)}.
 Then   for any $ (k, p, q)\in (1,\infty] \times [1,\infty)\times (0,\infty)$, there exists a constant $c\in (0,\infty)$ such that
\begin{equation}\label{UP1}
\big(\E^\nu[\W_p^q(\mu_t, \mu)] \big)^{\ff 1 q} \le c \|h_\nu\|_{L^k(\mu)}^{\ff 1 q} \gg_d(t),\ \ t\ge 2, \ \nu= h_\nu\mu\in \scr P \text{ with } h_\nu \in L^k(\mu), \end{equation}
where $\nu= h_\nu\mu$ means $\ff{\d\nu}{\d\mu}=h_\nu.$
 \end{thm}

\subsection{Lower bound estimate}

To derive sharp lower bound estimates, we make the following assumption {\bf (B)} which holds in particular for the (reflecting) diffusion operator  $\hat L:=\DD+\nn V $   on a $d$-dimensional compact connected Riemannian manifold, since in this case conditions \eqref{EG2} and \eqref{CT2} are well known,
and  the other conditions have been verified by   \cite[Lemma 5.2]{WZ}.  For $M$ being a smooth domain in $\R^d$,  \eqref{SA} is known as  Sard's lemma (see \cite[p130, Excercise 5.5]{Gry}) and
has been discussed in \cite[Section 3.1.6]{BLG}. The function $f_\xi$ in \eqref{LU} is called   Lusin's approximation  of $h$ (see \cite{AF, Liu}).

\beg{enumerate}\item[{\bf (B)}]    Let $\{\ll_i\}_{i\ge 0}$ be all eigenvalues of  $-\hat L$ listed in the increasing order with multiplicities. There exist  constants $\kk>0$ and $d\in [1,\infty)$ such that
\beq\label{EG2}  \ll_i \le \kk\,  i^{\ff  2 d},\ \ \ i\ge 0,\end{equation}
\beq\label{CT2} \W_1(\nu \hat P_t,\mu)\le \kk\, \W_1(\nu,\mu),\ \ t\in [0,1], \, \nu\in \scr P.\end{equation}
Moreover,   for any $f\in\D(\hat\EE)$,
\beq\label{SA} \mu\big(\{|\nn f|>0, f=0\}\big)=0,\end{equation}
and there exists a constant $c>0$ independent of $f$ such that
\beq\label{LU}
\mu \big(   f  \ne f_\xi   \big) \leq \frac{c}{\xi^2} \int_M |\nabla f|^2 \d \mu,\ \ \xi>0
\end{equation} holds for   a family of functions $\{f_\xi:\xi>0\}$ on $M$ with   $\|\nn f_\xi \|_{\infty} \leq \xi$.
\end{enumerate}

 \beg{thm}[\cite{WWZ}]Assume {\bf (A)} and {\bf (B)}. Then for any  $q\in (0,\infty)$, we have for large $t>0$,
 $$\inf_{\nu\in \scr P} \E^\nu[\W_1^q(\mu_t,\mu)]   \succeq \gg_d(t)^q.$$
 \end{thm}
 
 Finally,  we present a lower bound estimate which also applies to infinite-dimensions and generalizes \cite[Proposition 4.2]{RE1}  for the finite-dimensional setting.

\beg{thm}\label{LB101} Let $\mu\in \scr P(E) $ such that
\beq\label{LB1} \sup_{x\in E}\mu(B(x,r)) \le \psi(r),\ \ r\ge 0 \end{equation} holds for some increasing function $\psi,$ where
$B(x,r):=\{y\in E: \rr(x,y)<r\}$. Then
for any $N\ge 1$ and any probability measure $\mu_N$ supported on a set of $N$ points in $E$,
\beq\label{LB2} \W_p(\mu_N,\mu) \ge 2^{-\ff 1 p}  \psi^{-1}\Big(\ff 1 {2N}\Big),\end{equation}
where $\psi^{-1}(s):=\sup\{r\ge 0: \psi(r)\le s\}, s\ge 0.$
\end{thm}
\beg{proof} Let $D={\rm supp}\mu_N$ which contains $N$ many points, so that from \eqref{LB1} we conclude that
 $D_r:= \cup_{x\in D} B(x,r)$ satisfies
$$\mu(D_r)\le \sum_{x\in D} \mu(B(x,r))\le N \psi(r),\  \ r\ge 0.$$
Therefore, for any $\pi\in \C(\mu_N, \mu)$, we get
$$ \int_{E\times E} \rr(x,y)^p \pi(\d x,\d y) \ge \int_{D\times D_r^c} r^p \pi(\d x,\d y) = r^p \mu(D_r^c) \ge r^p\{1-N\psi(r)\},\  \ r\ge 0.$$
Taking $r= \psi^{-1}(1/(2N))$ we derive 
$$\W_p(\mu,\nu)^p \ge \sup_{r\ge 0} r^p[1-N\psi(r)]\ge \ff 1 2 \big\{\psi^{-1}(1/(2N))\big\}^p.$$
\end{proof}

To apply Theorem \ref{LB101} to the empirical measure $\mu_t$, we only need to compare $\mu_t$ with the discretize empirical measure
$$\mu_{t,N}:=\ff 1 N\sum_{i=1}^N \dd_{X_{i t/N}},$$
with suitable choice of $N=N_t\to\infty$ as $t\to\infty$. 

\paragraph{Acknowledgement.}  The author would like to thanks the referee for helpful comments and corrections.

\end{document}